\documentclass[11pt,reqno]{amsart}
\usepackage{amsmath,amsfonts,amssymb,amsthm}
\usepackage{psfig}
\newlength{\figboxwidth}
\renewcommand{\setminus}{{\smallsetminus}}

\newcommand{\neigh}{\rm neigh} 

\def\@ifundefined#1#2#3%
  {\expandafter\ifx\csname#1\endcsname\relax#2\else#3\fi}

\@ifundefined{theoremstyle}{
}{
\theoremstyle{plain} 
}
\newtheorem{theorem}{Theorem}[section]

\newtheorem{proposition}[theorem]{Proposition}
\newtheorem{lemma}[theorem]{Lemma}

\newtheorem{corollary}[theorem]{Corollary}

\@ifundefined{theoremstyle}{
}{
\theoremstyle{definition} 
}

\newtheorem{remark}[theorem]{Remark}

\newlength{\halfbls}

\newcommand{\NN}{{\mathbb{N}}}
\newcommand{\cutoff}{{C}}

\catcode`\@=11
\long\def\@savemarbox#1#2{\global\setbox#1\vtop{\hsize\marginparwidth 
  \@parboxrestore\tiny\raggedright #2}}
\catcode`\@=12
\newcommand{\CC}{{\mathcal C}}

\begin{document}

\title{The pants complex has only one end}
\author{Howard Masur}
\author{Saul Schleimer}
\address{\hskip-\parindent
        Howard Masur\\
        Department of Mathematics, UIC\\
        851 South Morgan Street\\
        Chicago, Illinois 60607}
\email{masur@math.uic.edu}
\address{\hskip-\parindent
        Saul Schleimer\\
        Department of Mathematics, UIC\\
        851 South Morgan Street\\
        Chicago, Illinois 60607}
\email{saul@math.uic.edu}

\date{\today}

\maketitle

\section{Definitions and statement of the main theorem}

The purpose of this note is to prove the following theorem:

\vspace{2mm}
\noindent
{\bf Theorem~\ref{Thm:PantsComplexHasOnlyOneEnd}.}
{\em 
Let $S$ be a closed orientable surface with genus $g(S) \geq 3$. Then
the pants complex of $S$ has only one end.  In fact, there are
constants $K = K(S)$ and $M_3 = M_3(S)$ so that, if $R > M_3$, any
pair of pants decompositions $P$ and $Q$, at distance greater than
$KR$ from a basepoint, can be connected by a path which remains at
least distance $R$ from the basepoint.
}
\vspace{1mm}

Recall that a {\em pants decomposition} of $S$ consists of $3g(S) - 3$
disjoint essential non-parallel simple closed curves on $S$.  Each
component of the complement of the curves is a sphere with $3$ holes
or pair of pants. Then the {\em pants complex} $C_P(S)$ is the metric
graph whose vertices are pants decompositions of $S$, up to isotopy.
Two vertices $P, P'$ are connected by an edge if $P, P'$ differ by an
{\em elementary move}.  In an elementary move all curves in the pants
are fixed except for one curve $\alpha$.  Removing $\alpha$, the
component of the complement of the remaining curves that contains
$\alpha$ is either a punctured torus or a sphere with $4$ holes.  Then
$\alpha$ is replaced by a curve $\beta$ contained in this domain that
intersects $\alpha$ minimally; in the punctured torus case once, and
in the sphere case twice.  All edges of $C_P(S)$ are assigned length
$1$. We let $d(\cdot,\cdot)$ be the distance function in $C_P(S)$.
The pants complex $C_P(S)$ is known to be
connected~\cite{Hatcher:Thurston}.

Recall that a metric space $(X,d)$ has {\em one end} if for any basepoint
$O \in X$ and any radius $R$ the complement of $B_R = B_R(O)$, the
ball of radius $R$ centered at $O$, has only one unbounded
component.  It is easy to see that the definition does not depend on
the choice of the point $O$.  Clearly having one end is a
quasi-isometry invariant.  So, following Brock~\cite{Brock}, our
theorem implies:

\begin{corollary}
\label{Cor:WPHasOnlyOneEnd}
Teichm\"uller space, equipped with the Weil-Peterson metric, has only
one end.
\end{corollary}

Finally, recall that the {\em curve complex} $\CC(S)$ is the complex
whose $k$-simplices consist of $k+1$ distinct isotopy classes of
essential simple closed curves on $S$ that have disjoint
representatives on $S$.  Or, in the case of a once-punctured torus and
four-punctured sphere, $\CC(S)$ is the Farey graph.  From the metric
point of view we will only interested in the $1$-skeleton of $\CC(S)$.
Each edge is assigned length $1$. We let $d_S(\cdot,\cdot)$ denote the
distance function in $\CC(S)$.

\section{The set of handle curves is connected}
\label{Sec:H(S)Connected}

In this section we prove two combinatorial facts.  First, the set of
handle curves in the curve complex is connected and second, any pants
decomposition is a bounded distance (in the pants complex) from a
decomposition containing a handle curve.

Again assume $S$ is a closed orientable surface with genus three or
greater.  We will call $\alpha$ a {\em handle curve} in $S$ if
$\alpha$ separates $S$ into two surfaces: the handle $S(\alpha) \cong
{\mathbb{T}^2 \setminus {\rm \{pt\}}}$ and the rest of the surface.

We will need the following result. It was first proved by Farb and
Ivanov~\cite{Farb:Ivanov} by different methods.  Another proof was
given by McCarthy and Vautaw~\cite{McCarthy:Vautaw} by methods similar
to ours.  We include a proof for completeness.

\begin{proposition}
\label{Prop:H(S)IsConnected}
If $g(s) \geq 3$, the subcomplex $h(S) \subset \CC(S)$ of handle
curves is connected.
\end{proposition}

\begin{remark}
Note that the hypothesis $g(S) \geq 3$ cannot be removed; it is easy
to check that the handle curves for a closed surface of genus $2$ do
not form a connected set.
\end{remark}

\begin{remark}
Note that Proposition~\ref{Prop:H(S)IsConnected} immediately implies that
the set of separating curves in $\CC(S)$ is also connected.
\end{remark}

\begin{remark}
We note that Proposition~\ref{Prop:H(S)IsConnected} may be
generalized to the case $\partial S \neq \emptyset$.  A still open
question is the higher connectivity of $h(S)$.
\end{remark}

Before we begin the proof we will require a bit of terminology.  Let
$i(\cdot,\cdot)$ denote the geometric intersection number of two
essential simple closed curves in $S$.  Also, if $\delta$ is a
dividing curve in $S$ we say that an arc $\beta'$ is a {\em wave} for
$\delta$ if $\beta' \cap \delta = \partial \beta'$ and $\beta'$ is
essential as a properly embedded arc in $S \setminus \delta$.  We say
that two waves $\beta'$ and $\beta''$ for $\delta$ {\em link} if
$\beta' \cap \beta'' = \emptyset$, both $\beta'$ and $\beta''$ meet
the same side of $\delta$, and $\partial \beta'$ separates $\partial
\beta''$ inside $\delta$.

Finally we define {\em double surgery} as follows.  Suppose we are
given a linking pair of waves $\beta'$ and $\beta''$ for an essential
dividing curve $\delta_0$.  Form the closed regular neighborhood $U =
\neigh(\delta_0 \cup \beta' \cup \beta'')$.  Let $\delta_1$ be the
component of $\partial U$ which is not homotopic in $S$ to $\delta_0$.
We say that $\delta_1$ is obtained from $\delta_0$ via double surgery
along $\beta'$ and $\beta''$.
Note that $\delta_1$ is necessarily a dividing curve and is disjoint
from $\delta_0$.  If the component of $S \setminus \delta_0$
containing $\beta ' \cup \beta''$ is not a handle then $\delta_1$ is
also essential.

We are now equipped to prove the proposition:

\begin{proof}[Proof of Proposition~\ref{Prop:H(S)IsConnected}]
Let $\alpha, \beta \in h(S)$ be handle curves and $S$ a closed
orientable surface of genus at least three.  Suppose that $\alpha$ and
$\beta$ are {\em tight}: $\alpha$ has been isotoped to make $|\alpha
\cap \beta| = i(\alpha,\beta)$.  If $i(\alpha,\beta) = 0$ then there
is nothing to prove.  If $i(\alpha,\beta)>0$ we will find a curve
$\gamma \in h(S)$ with $i(\gamma ,\alpha) = 0$ and
$i(\gamma,\beta)<i(\alpha,\beta)$.  By induction, $\gamma$ will be
connected to $\beta$ in $h(S)$, proving the proposition.

We find $\gamma$ via the following inductive procedure.  Recall that
$S(\alpha)$ is the handle which $\alpha$ bounds.  To begin, we define
$\delta_0 \subset S \setminus S(\alpha)$ to be a parallel copy of
$\alpha$, still intersecting $\beta$ tightly.  At stage $k$ by
induction we will be given an essential dividing curve $\delta_k$
where
\begin{itemize}
\item $i(\alpha,\delta_k) = 0$,
\item $\delta_k$ is tight with respect to $\beta$, and
\item $i(\delta_k,\beta)<i(\delta_{k-1},\beta)$, if $k > 0$. 
\end{itemize}
Let $T_k$ be the component of $S \setminus \delta_k$ which does not
contain $\alpha$.  If $T_k$ is a handle, then we take $\gamma =
\delta_k$ and we are done with the inductive procedure.  If
$i(\delta_k,\beta) = 0$ then we may take $\gamma$ to be any handle
curve inside $T_k$.  As this $\gamma$ satisfies $i(\alpha, \gamma) =
i(\beta,\gamma) = 0$ this would finish the proposition.  From now on
we assume that $T_k$ is not a handle and that $i(\delta_k, \beta) >
0$.

We now attempt to do a double surgery of $\delta_k$ into $T_k$.
Either we will find $\gamma$ directly or the curve resulting from
double surgery, $\delta_{k+1}$, will satisfy the induction hypothesis.
As the geometric intersection with $\beta$ is always decreasing, this
procedure will stop after finitely many steps yielding the desired
handle curve.

So all that remains is to do the double surgery.  Recall that we are
given $\alpha, \beta$ tightly intersecting handle curves and we are
also given $\delta_k$ satisfying the induction hypotheses.  Recall also
that $T_k$ is the component of $S \setminus \delta_k$ which does not
contain $\alpha$.  Recall $T_k$ is not a handle and that
$i(\delta_k,\beta) > 0$.

Suppose further that $\beta', \beta'' \subset \beta \cap T_k$ are
linking waves for $\delta_k$.  As described above we may form
$\delta_{k+1}$ via a double surgery along $\beta'$ and $\beta''$.
Isotope $\delta_{k+1}$, in the complement of $\delta_k$, to be tight
with respect to $\beta$.  As noted in the definition of double
surgery, $\delta_{k+1}$ is an essential
dividing curve
which is disjoint from $\alpha$.
Finally note that $i(\delta_{k+1},\beta)\leq i(\delta_k,\beta) - 4$.
Thus all of the induction hypotheses are satisfied.

Suppose now that we cannot find linking waves among the arcs $\beta
\cap T_k$.  Choose instead an {\em outermost} wave $\beta' \subset
\beta \cap T_k$: that is, there exists an arc $\delta_k' \subset
\delta_k$ such that $\delta_k' \cap \beta = \partial \delta_k' =
\partial \beta'$.


Here there are two remaining cases.  If $\delta_k' \cup \beta'$ is a
separating curve take $\delta_{k+1} = \delta_k' \cup \beta'$ and note
that the induction hypotheses are easily verified.
The final possibility is that $\delta_k' \cup \beta'$ is not
separating.
In this case choose a properly embedded essential arc $\beta'' \subset
T_k$ such that $\beta'' \cap \beta = \emptyset$ and $|\beta'' \cap
\delta_k'| = 1$.  Then $\beta'$ and $\beta''$ link.  Do a double
surgery along these waves to obtain $\delta_{k+1}$.  Isotope
$\delta_{k+1}$, in the complement of $\delta_k$, to be tight with
respect to $\beta$.  Again, all of the induction hypotheses are easily
verified, as we have $i(\delta_{k+1},\beta) \leq i(\delta_k,\beta)
- 2$.  This completes the second induction step and hence completes
the proof of Proposition~\ref{Prop:H(S)IsConnected}.
\end{proof}

We also require

\begin{lemma}
There is a constant $M_3 = M_3(S)$ such that the pants decompositions
containing a handle curve are $M_3$-dense in the space of all pants
decompositions
\end{lemma}

\begin{proof}
The mapping class group acts co-compactly on the space of pants
decompositions.
\end{proof}

\section{Subsurface projections and distances}
\label{Sec:Projections}

Here we give two lemmas studying the pants complex.  The first gives a
condition for a pants decomposition to lie outside of a large ball
about the origin in $C_P(S)$ while the second provides us with useful
paths laying outside of such a ball.  This uses techniques developed
by Masur and Minsky~\cite{Masur:Minsky}.

Given a subsurface $W\subset S$ and a curve $\gamma$ that intersects
$W$, we may define a projection of $\gamma$ to $\CC(W)$ which
associates to $\gamma$ a collection of curves in $\CC(W)$.  Namely,
the intersections of $\gamma$ with $W$ fall into finitely many
homotopy classes of disjoint arcs (and curves) relative to $\partial
W$.  Let $\alpha$ be any such arc or curve.  Let $U$ be a regular
neighborhood of $\alpha \cup \partial W$. Then $\partial U \setminus
\partial W$ is a curve in $\CC(W)$ if $\alpha$ is a curve or is an arc
connecting distinct components of $\partial W$.  Otherwise $\partial U
\setminus \partial W$ is a pair of curves in $\CC(W)$.  Clearly this
curve or pair of curves in $\CC(W)$ depend only on the homotopy class
of $\alpha$.  If $\alpha,\beta$ are disjoint homotopy classes of
arcs, then any pair of curves $\alpha',\beta'$ built out of this
surgery satisfy $d_W(\alpha',\beta')\leq 2$ (\cite{Masur:Minsky}
Lemma 2.3).  Thus we can define $\pi_W (\gamma)$ to be the
corresponding subset of diameter at most $2$ in $\CC(W)$.

Similarly, given a pants decomposition $P$ we may project each curve
of $P$ that intersects $W$ into $W$.  We denote the resulting image
set in $\CC(W)$, which has diameter at most $4$, by $\pi_W (P)$.  By
$d_W(P, P')$ we mean the distance in the curve complex of $W$
between the sets $\pi_W(P)$ and $\pi_W(P')$.

Let $[x]_\cutoff$ be the function on $\NN$ giving zero if $x <
\cutoff$ and giving $x$ if $x \geq \cutoff$.  We will need the
following result from~\cite{Masur:Minsky} (see Theorem~6.12 and
Section~8 of that paper):

There is a constant $C_0=C_0(S) \geq 1$ such that for any $\cutoff
\geq C_0$ there are constants $M_1 = M_1(\cutoff) \geq \cutoff$ and
$M_2 = M_2(\cutoff) \geq 0$ with the following property: for any pants
decompositions $P, P'$ we have
\begin{equation}
\label{eq:projections}
\frac{1}{M_1} \sum_V [d_V(P, P')]_\cutoff - M_2
\leq d(P, P')
\leq M_1 \sum_V [d_V(P, P')]_\cutoff + M_2,
\end{equation}
where the sums range over subsurfaces $V \subset S$ with essential
boundary and where $V$ is not an annulus nor is it a thrice-punctured
sphere.

Fix now such a $\cutoff > 2$.  It follows from
equation~(\ref{eq:projections}) that the projections that are at the
critical values $\cutoff, \cutoff+1$ cannot account for the entire
distance. Namely there are constants $c = c(\cutoff) > 1$ and
$M_4 = M_4(\cutoff) > 0$ such that
\begin{equation}
\label{eq:critical}
\sum_V [d_V(P, P')]_\cutoff \leq 
                c \cdot \sum_V [d_V(P, P')]_{\cutoff + 2} + M_4.
\end{equation}

Choose $K = K(\cutoff) > 0$ so that for all $R \geq 1$,
\begin{equation}
\label{eq:choose}
\frac{1}{2cM_1} \big( (K-1)R-M_2- M_1 M_4-cM_1^2(R+M_2) \big)
              > M_1(R+M_2)
\end{equation}
Also, choose a basepoint $O \in C_P(S)$ and let $B_R = B_R(O)$ be 
the ball of radius $R$ centered at $O$.  

\begin{lemma}
\label{Lem:Points}
Fix a handle curve $\alpha$ and some curve $\alpha'' \subset
S(\alpha)$ satisfying $d_{S(\alpha)}(O, \alpha'') > M_1(R + M_2)$.
For any pants decomposition $P$ containing $\alpha''$ we have $P \notin
B_R$.
\end{lemma}

\begin{proof}
Note that $d_{S(\alpha)}(O, \alpha'') \geq C$. So, by the left
inequality of (\ref{eq:projections}) any pants $P$ containing
$\alpha''$ has
\begin{displaymath}
d(P, O) \geq 
\frac{1}{M_1} [d_{S(\alpha)}(P, O)]_C - M_2 \geq
\frac{1}{M_1} d_{S(\alpha)}(\alpha'', O) - M_2 > R.
\end{displaymath}
\end{proof}

As the Farey graph for $S(\alpha)$ has infinite diameter, and as the
diameter of $\pi_{S(\alpha)}(O)$ is bounded, such curves $\alpha''$
exist in abundance.  We now turn to the existence of paths lying
outside of the $R$-ball about the basepoint.

\begin{lemma}
\label{Lem:Paths}
Suppose $P_0$ is a pants decomposition of $S$ such that $P_0 \notin
B_{(K-1)R}$ and $P_0$ contains a curve $\alpha$ which bounds a handle
$S(\alpha)$. Then there is a path $P_t$ starting at $P_0$ such that
\begin{itemize}
\item  for all $t$, $P_t|(S \setminus S(\alpha))=P_0|(S\setminus S(\alpha))$.
\item for all $t$, $P_t \notin B_R$
\item The endpoint of the path, $P_1$, contains a curve $\alpha''
\subset S(\alpha)$ which does not appear in any pants decomposition in
$B_R$.
\end{itemize}
\end{lemma}

\begin{proof}
Let $\alpha' \in P_0$ be the curve strictly contained in $S(\alpha)$.
Consider a geodesic segment in the Farey graph connecting $\alpha'$ to
$\beta \in \pi_{S(\alpha)}(O)$, where $\beta$ is chosen as close as
possible to $\alpha'$.  Extend this segment through $\alpha'$ to a
geodesic ray $L$ in the direction opposite $\beta$.  The ray $L$ meets
the segment only at $\alpha'$.  Move along $L$ distance more than
$M_1(R + M_2)$ from $\alpha'$ to a point $\alpha''$.  Let $P_t$ be the
path obtained by making elementary moves along the curves in $L$ and
fixing the pants in $S \setminus S(\alpha)$.

Suppose first that $d_{S(\alpha)}(\beta, \alpha') > M_1(R + M_2)$.  Then by
Lemma~\ref{Lem:Points} any pants $\hat P$ containing any $\alpha_t \in
L$ has $d(\hat P, O) > R$.  So Lemma~\ref{Lem:Paths} holds in this case.

Next suppose that $d_{S(\alpha)}(\beta, \alpha')\leq M_1(R+M_2)$.
Then by (\ref{eq:projections}) and (\ref{eq:critical})

$$(K - 1) R \leq d(P_0, O) \leq M_1 \sum_V [d_V(P_0, O)]_C + M_2 \leq$$

$$\leq c M_1 \sum_V [d_V(P_0,0)]_{\cutoff+2} + M_1 M_4 + M_2 \leq$$

$$\leq c M_1 \sum_{V\neq S(\alpha)} [d_V(P_0,O)]_{\cutoff+2} 
                             + c M_1^2 (R + M_2) + M_1 M_4 + M_2.$$ 

Let $V$ be any subsurface disjoint from $S(\alpha)$.  Since $P_t$ is
constant in $V$, the projection $\pi_V(P_t)$ is constant.  Now let $V$
be a subsurface that intersects $S(\alpha)$ or strictly contains
$S(\alpha)$.  Since $\alpha \in P_t$, it follows that $\pi_V(P_t)$
contains $\pi_V(\alpha)$.  Since each $\pi_V(P_t)$ has diameter at
most $2$, $d_V(P_t,O) \geq d_V(P_0,O) - 2$.  Thus for any subsurface
$V$ not isotopic to $S(\alpha)$, as $\cutoff > 2$, we have
$[d_V(P_t,O)]_\cutoff \geq \frac{1}{2} [d_V(P_0,O)]_{\cutoff + 2}$.
Thus,

\begin{displaymath}
\sum_{V\neq S(\alpha)} [d_V(P_t,O)]_\cutoff \geq 
\frac{1}{2} \sum_{V\neq S(\alpha)} [d_V(P_0,0)]_{\cutoff+2} \geq
\end{displaymath}
\begin{displaymath} 
\geq \frac{1}{2cM_1}\big((K-1)R-M_2-M_1M_4-cM_1^2(R+M_2)\big) 
  > M_1(R+M_2),
\end{displaymath}
the last inequality following from (\ref{eq:choose}).  So, by
(\ref{eq:projections}),
\begin{displaymath}
d(P_t,O) > R.
\end{displaymath}  

Finally, as the pants decomposition $P_1$ contains $\alpha''$ and
$$d_{S(\alpha)}(O, \alpha'') \geq d_{S(\alpha)}(\alpha', \alpha'')
> M_1(R + M_2)$$ we have $P_1 \notin B_R$, by
Lemma~~\ref{Lem:Points}.
\end{proof}

\section{Proof of the theorem}

Recall the statement:

\begin{theorem}
\label{Thm:PantsComplexHasOnlyOneEnd}
Let $S$ be a closed orientable surface with genus $g(S) \geq 3$. Then
the pants complex of $S$ has only one end.  In fact, there are
constants $K = K(S)$ and $M_3 = M_3(S)$ so that, if $R > M_3$, any
pair of pants decompositions $P$ and $Q$, at distance greater than
$KR$ from a basepoint, can be connected by a path which remains at
least distance $R$ from the basepoint.
\end{theorem}

\begin{proof}
We take $M_3$ as defined in Section~\ref{Sec:H(S)Connected} and $K$
as defined in Section~\ref{Sec:Projections}.

First move $P$ and $Q$ a distance at most $M_3 < R$ to obtain pants
decompositions $P_0$ and $Q_0$ which contain handle curves $\alpha_P
\in P_0$, $\alpha_Q \in Q_0$ and such that $P_0, Q_0 \notin B_{(K-1)R}$.

Apply Lemma~\ref{Lem:Paths} twice in order to connect $P_0$ and $Q_0$
to pants decompositions $P_1$ and $Q_1$ satisfying all of the
conclusions of the lemma.  Let $\alpha_P''$ and $\alpha_Q''$ be the
curves lying in the handles $S(\alpha_P) \subset P_1$ and $S(\alpha_Q)
\subset Q_1$ respectively.

We must now construct a path from $P_1$ to $Q_1$.  Consider first the
case where $\alpha_P \neq \alpha_Q$.  

Applying Proposition~\ref{Prop:H(S)IsConnected} we connect $\alpha_P
\in P_1$ and $\alpha_Q \in Q_1$ by a path $\{\alpha_i\}_{i = 1}^n$ of
handle curves in $\CC(S)$, the curve complex of $S$.  Here we have
$\alpha_1 = \alpha_P$, $\alpha_n = \alpha_Q$, and $n > 1$.  Note that
in this step the hypothesis $g(S) > 2$ is used.  Choose, for $i \in
\{2, 3, \ldots n - 1\}$, any curve $\alpha_i'' \subset S(\alpha_i)$
such that $d_S(\alpha_i) (O, \alpha_i'') > M_1(R + M_2)$.  Set
$\alpha_1'' = \alpha_P''$ and $\alpha_n'' = \alpha_Q''$.  Let $P_n =
Q_1$.

Inductively, we connect $P_i$ by a path to $P_{i+1}$ where, first,
every pants decomposition in the path contains $\alpha_i$ and
$\alpha_i''$ and, second, $P_{i+1}$ also contains $\alpha_{i+1}$ and
$\alpha_{i+1}''$.  (This is possible because $C_P(S \setminus
S(\alpha_i))$ is connected.)  By Lemma~\ref{Lem:Points} this path lies
outside of the ball of radius $R$ and we are done.

In the case which remains we have $\alpha_P = \alpha_Q$.  So there is
no need for Proposition~\ref{Prop:H(S)IsConnected}.  Instead we choose
any handle curve $\beta$ which is disjoint from $\alpha_P$.  Note that
$\beta$ exists as $g(S) > 2$.  Choose also any $\beta''$ satisfying
the hypothesis of Lemma~\ref{Lem:Points}.  We now consider the
sequence $\alpha_P, \beta, \alpha_Q$ as a path of length two in the
curve complex and connect $P_1$ to $Q_1$ as in the previous paragraph.
This completes the proof.
\end{proof}

\end{document}